\documentstyle[11pt]{article}

\newtheorem{eg}{{\bf Example}}
\newtheorem{lemma}{{\bf Lemma}}[section]
\newtheorem{theo}{{\bf Theorem}}
\newtheorem{prop}{{\bf Proposition}}[section]

\newtheorem{remark}{{\bf Remark}}

\font\bbb=msbm10 scaled\magstep1

\newcommand{\ZZ}{\mbox{\bbb Z}}

\def\scoll{\mbox{{\scriptsize
$\,\searrow$}$\!\!\stackrel{{}^{\rm s}}{}\,\,$}}

\def\coll{\mbox{\scriptsize $\,\,\searrow\,$}}

\newcommand{\La}{\Lambda}

\newcommand{\si}{\sigma}

\begin{document}

\title{\bf Uniqueness of Walkup's 9-vertex 3-dimensional Klein bottle}
\author{{\bf Bhaskar Bagchi}$^{\rm a}$ and {\bf Basudeb Datta}$^{\rm b}$
}

\date{}

\maketitle

\vspace{-5mm}

\noindent {\small $^{\rm a}$ Theoretical Statistics and
Mathematics Unit, Indian Statistical Institute,  Bangalore
560\,059, India.

\smallskip

\noindent $^{\rm b}$ Department of Mathematics, Indian Institute
of Science, Bangalore 560\,012,  India.}\footnotetext{{\em E-mail
addresses:} bbagchi@isibang.ac.in (B. Bagchi),
dattab@math.iisc.ernet.in (B. Datta).}

\begin{center}


\end{center}

\hrule

\bigskip

 {\small

\noindent {\bf Abstract}

\medskip

Via a computer search, Altshuler and Steinberg found that there
are $1296 +1$ combinatorial 3-manifolds on nine vertices, of
which only one is non-sphere. This exceptional 3-manifold $K^3_9$
triangulates the twisted $S^{\,2}$-bundle over $S^{\,1}$. It was
first constructed by Walkup. In this paper, we present a
computer-free proof of the uniqueness of this non-sphere
combinatorial 3-manifold. As opposed to the computer-generated
proof, ours does not require wading through all the 9-vertex
3-spheres. As a preliminary result, we also show that any
9-vertex combinatorial 3-manifold is equivalent by proper
bistellar moves to a 9-vertex neighbourly 3-manifold. }

\bigskip

{\small

 \noindent 2000 Mathematics Subject Classification. 57Q15,
57R05.

\smallskip

\noindent Keywords. Combinatorial 3-manifolds,  pl manifolds,
Bistellar moves.

}

\bigskip

\hrule

\section{Introduction and results}

Recall that a {\em simplicial complex} is a collection of
non-empty finite sets (sets of {\em vertices}) such that every
non-empty subset of an element is also an element.  For $i \geq
0$, the elements of size $i+1$ are called the {\em $i$-simplices}
(or {\em $i$-faces}) of the complex. For a simplicial complex
$K$, the maximum of $k$ such that $K$ has a $k$-simplex is called
the {\em dimension} of $K$ and is denoted by $\dim(K)$. If any
set of $\lfloor\frac{d}{2}\rfloor +1$ vertices form a face of a
$d$-dimensional simplicial complex $K$, then one says that $K$ is
{\em neighbourly}.

All the simplicial complexes considered here are finite. The
vertex-set of a simplicial complex $K$ is denoted by $V(K)$. If
$K$, $L$ are two simplicial complexes, then a {\em simplicial
isomorphism} from $K$ to $L$ is a bijection $\pi : V(K) \to V(L)$
such that for $\si\subseteq V(K)$, $\si$ is a face of $K$ if and
only if $\pi (\si)$ is a face of $L$. Two complexes $K$, $L$ are
called {\em isomorphic} when such an isomorphism exists. We
identify two simplicial complexes if they are isomorphic.

A simplicial complex is usually thought of as a prescription for
construction of a topological space by pasting geometric
simplices. The space thus obtained from a simplicial complex $K$
is called the {\em geometric carrier} of $K$ and is denoted by
$|K|$. If a topological space $X$ is homeomorphic to $|K|$ then
we say that $K$ is a {\em triangulation} of $X$. A {\em
combinatorial $d$-manifold} is a triangulation of a closed pl
$d$-manifold (see Section 2 for more).

For a set $V$ with $d + 2$ elements, let $S$ be the simplicial
complex whose faces are all the non-empty proper subsets of $V$.
Then $S$ triangulates the $d$-sphere. This complex is called the
{\em standard $d$-sphere} and is denoted by $S^{\,d}_{d + 2}(V)$
or simply by $S^{\,d}_{d + 2}$.

If $\si$ is a face of a simplicial complex  $K$ then the {\em
link} of $\si$ in $K$, denoted by ${\rm lk}_K(\sigma)$ (or simply
by ${\rm Lk}(\sigma)$), is by definition the simplicial complex
whose faces are the faces $\tau$ of $K$ such that $\tau$ is
disjoint from $\si$ and $\si\cup\tau$ is a face of $K$.

Let $M$ be a $d$-dimensional simplicial complex. If $\alpha$ is a
$(d - i)$-face of $M$, $0 < i \leq d $, such that ${\rm
lk}_M(\alpha) = S^{i - 1}_{i + 1}(\beta)$ and $\beta$ is not a
face of $M$ (such a face $\alpha$ is said to be a {\em removable}
face of $M$) then consider the complex (denoted by
$\kappa_{\alpha}(M)$) whose set of maximal faces is $\{\sigma :
\sigma \mbox{ a maximal face of } M, \alpha \not\subseteq \sigma\}
\cup \{ \beta \cup \alpha \setminus \{v\} : v \in \alpha\}$. The
operation $\kappa_{\alpha} : M \mapsto \kappa_{\alpha}(M)$ is
called a {\em bistellar $i$-move}. For $0< i < d$, a bistellar
$i$-move is called a {\em proper bistellar move}. In \cite{as},
Altshuler and Steinberg found from their computer search that all
the 9-vertex 3-spheres are equivalent via a finite sequence of
proper bistellar moves. Here we prove:

\begin{theo}$\!\!${\bf .} \label{t1}
Every $9$-vertex combinatorial $3$-manifold is obtained from a
neighbourly $9$-vertex combinatorial $3$-manifold by a sequence
of $($at most $10)$ bistellar $2$-moves.
\end{theo}

In \cite{a1}, Altshuler has shown that every combinatorial
3-manifold with at most $8$ vertices is a combinatorial 3-sphere.
(This is also a special case of a more general result of Brehm
and K\"{u}hnel \cite{bk}.) Via a computer search, Altshuler and
Steinberg found (in \cite{as}) that there are $1297$
combinatorial 3-manifolds on nine vertices, of which only one
(namely, $K^3_9$ of Example \ref{e1} below) is non-sphere. Here,
we present a computer-free proof of this fact. More explicitly, we
prove\,:

\begin{theo}$\!\!${\bf .} \label{t2}
Up to simplicial isomorphism, there is a unique $9$-vertex
non-sphere combinatorial $3$-manifold, namely $K^3_9$.
\end{theo}

Note that Theorem \ref{t2} was a key ingredient in the proof of
the main result (viz, non-existence of complementary
pseudomanifolds of dimension 6) in \cite{bd4}. The proof of
Theorem \ref{t2} presented here makes the result of \cite{bd4}
totally independent of machine computations. This was one of the
prime motivations for the present paper.

\section{Preliminaries}

For $i = 1, 2, 3$, the $i$-faces of a simplicial complex $K$ are
also called the {\em edges}, {\em triangles} and {\em tetrahedra}
of $K$, respectively. A simplicial complex $K$ is called {\em
connected} if $|K|$ is connected. For a simplicial complex $K$,
if $U \subseteq V(K)$ then $K[U]$ denotes the induced subcomplex
of $K$ on the vertex-set $U$. If the number of $i$-simplices of a
$d$-dimensional simplicial complex  $K$ is $f_i(K)$, then the
vector $f= (f_0, \dots, f_d)$ is called the {\em $f$-vector} of
$K$ and the number $\chi(K) := \sum_{i= 0}^{d}(- 1)^i f_i(K)$ is
called the {\em Euler characteristic} of $K$.

For a face $\sigma$ in a simplicial complex $K$, the number of
vertices in ${\rm lk}_K(\sigma)$ is called the {\em degree} of
$\sigma$ in $K$ and is denoted by $\deg_K(\sigma)$. The induced
subcomplex $C(\sigma, K)$ on the vertex-set $V(K) \setminus
\sigma$ is called the {\em simplicial complement} of $\sigma$ in
$K$.

By a {\em subdivision} of a simplicial complex  $K$ we mean a
simplicial complex $K^{\,\prime}$ together with a homeomorphism
from $|K^{\,\prime}|$ onto $|K|$ which is facewise linear. Two
complexes $K$ and $L$ are called {\em combinatorially equivalent}
(denoted by {\em $K \approx L$}) if they have isomorphic
subdivisions. So, $K \approx L$ if and only if $|K|$ and $|L|$
are pl homeomorphic. If a simplicial complex $X$ is
combinatorially equivalent to $S^{\,d}_{d + 2}$ then it is called
a {\em combinatorial $d$-sphere}. A simplicial complex $K$ is
called a {\em combinatorial $d$-manifold} if the link of each
vertex is a combinatorial $(d - 1)$-sphere. Thus, a simplicial
complex $K$ is a combinatorial $d$-manifold if and only if $|K|$
is a closed pl $d$-manifold with the pl structure induced from
$K$ (see \cite{rs}).

A {\em graph} is an 1-dimensional simplicial complex. For $n\geq
3$, the $n$-vertex combinatorial 1-sphere ({\em $n$-cycle}) is
the unique $n$-vertex 1-pseudomanifold and is denoted by
$S^{\,1}_n$. A {\em coclique} in a graph is a set of pairwise
non-adjacent vertices.

A simplicial complex $K$ is called {\em pure} if all the maximal
faces of $K$ have the same dimension.  A maximal face in a pure
simplicial complex is also called a {\em facet}. For a pure
$d$-dimensional simplicial complex $K$, let $\La (K)$ be the graph
whose vertices are the facets of $K$, two such vertices being
adjacent in $\La (K)$ if and only if the corresponding facets
intersect in a $(d - 1)$-simplex. A $d$-dimensional pure
simplicial complex $K$ is called a {\em $d$-pseudomanifold} if
each $(d-1)$-face is contained in exactly two facets of $K$ and
$\La(K)$ is connected. If the link of a vertex in a
pseudomanifold is not a combinatorial sphere then it is called a
{\em singular vertex}. Clearly, any connected combinatorial
$d$-manifold is a $d$-pseudomanifold without singular vertices.
Since a connected $(d+1)$-regular graph has no $(d+1)$-regular
proper subgraph, a $d$-pseudomanifold has no proper
$d$-dimensional sub-pseudomanifold.

For two simplicial complexes $K$, $L$ with disjoint vertex-sets,
the {\em join} $K \ast L$ is the simplicial complex $K\cup L\cup
\{\sigma \cup \tau ~ : ~ \sigma \in K, \tau\in L\}$. Clearly, if
both $K$ and $L$ are pseudomanifolds then $K\ast L$ is a
pseudomanifold of dimension $\dim(K) + \dim(L) +1$.

Let $K$ be an $n$-vertex $d$-pseudomanifold. If $u$ is a vertex
of $K$ and $v$ is not a vertex of $K$ then consider the pure
simplicial complex $\Sigma_{uv}K$ on the vertex set $V(K) \cup
\{v\}$ whose set of facets is $\{\sigma \cup \{u\} : \sigma$ is a
facet of $K$ and $u \not\in \sigma\} \cup \{\tau \cup \{v\} :
\tau$ is a facet of $K\}$. Then $\Sigma_{uv}K$ is a
$(d+1)$-pseudomanifold and is called the {\em one-point
suspension} of $K$ (see \cite{bd2}). It is easy to see that the
links of $u$ and $v$ in $\Sigma_{uv}K$ are isomorphic to $K$.

\begin{eg}$\!\!${\bf .} \label{e1}
{\rm For $d \geq 2$, let $K^{\,d}_{2d+3}$ be the $d$-dimensional
pure simplicial complex whose vertices are the vertices of the
$(2d+3)$-cycle $S^{\,1}_{2d+3}$ and the facets are the sets of
$d+1$ vertices obtained by deleting an interior vertex from the
$(d+2)$-paths in the cycle. The simplicial complex
$K^{\,d}_{2d+3}$ is a combinatorial $d$-manifold. Indeed, it was
shown in \cite{k1} that $K^{\,d}_{2d+3}$ triangulates
$S^{\,d-1}\times S^{\,1}$ for $d$ even, and it triangulates the
twisted product of $S^{\,d-1}$ and $S^{\,1}$ for $d$ odd. In
particular, $K^{\,3}_9$ triangulates the twisted product of
$S^{\,2}$ and $S^{\,1}$ (often called the {\em $3$-dimensional
Klein bottle}). It was first constructed by Walkup in \cite{w}. }
\end{eg}

\begin{eg}$\!\!${\bf .} \label{e2}
{\rm Some combinatorial $2$-spheres on $5$, $6$ and $7$ vertices.}
\end{eg}

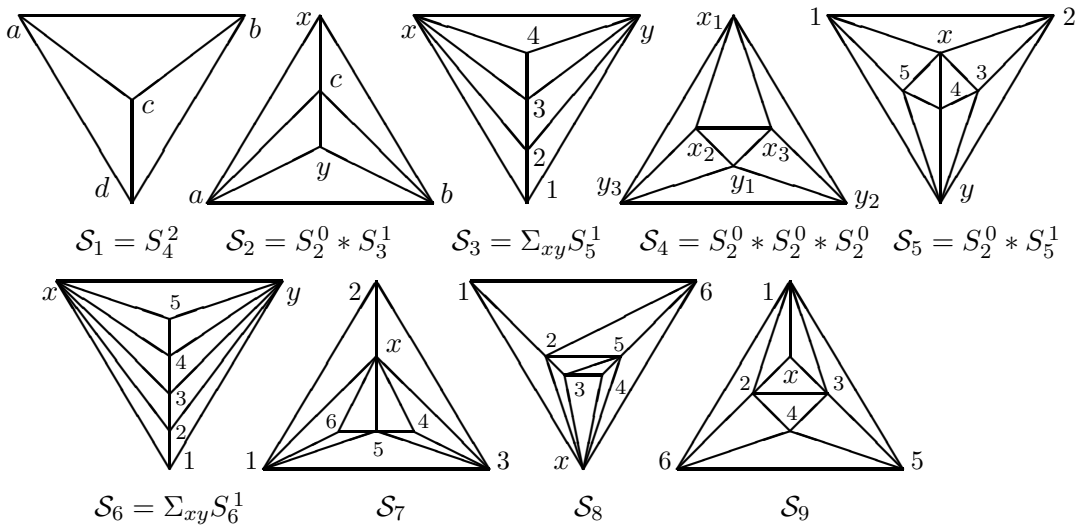
\setlength{\unitlength}{2.5mm}

\begin{picture}(55,13.5)(-9,0)

\thicklines

\put(-4,3){\line(-3,5){6}} \put(-4,3){\line(0,1){5.5}}
\put(-4,3){\line(3,5){6}} \put(-10,13){\line(1,0){12}}

\put(-10,13){\line(4,-3){6}} \put(2,13){\line(-4,-3){6}}

\put(-10.7,11.8){$a$} \put(2.2,11.8){$b$} \put(-3.5,7.8){$c$}
\put(-6,3.3){$d$}

\put(-7,0.5){\mbox{${\cal S}_1= S^{\,2}_4$}}


\put(0,3){\line(3,5){6}} \put(0,3){\line(1,1){6}}
\put(0,3){\line(2,1){6}} \put(0,3){\line(1,0){12}}
\put(6,6){\line(0,1){7}}

\put(12,3){\line(-3,5){6}} \put(12,3){\line(-1,1){6}}
\put(12,3){\line(-2,1){6}}

\put(-1,3){$a$} \put(12.4,3){$b$} \put(6.5,9.2){$c$}
\put(4.7,12.5){$x$} \put(5.8,4.5){$y$}

\put(1,0.5){\mbox{${\cal S}_2 = S^{\,0}_2\ast S^{\,1}_3$}}


\put(17,3){\line(-3,5){6}} \put(17,3){\line(0,1){8}}
\put(17,3){\line(3,5){6}} \put(11,13){\line(1,0){12}}

\put(11,13){\line(5,-6){6}} \put(23,13){\line(-5,-6){6}}
\put(11,13){\line(4,-3){6}} \put(23,13){\line(-4,-3){6}}
\put(17,11){\line(-3,1){6}} \put(17,11){\line(3,1){6}}

\put(10.3,11.8){$x$} \put(23,11.8){$y$} \put(17,11.5){\small 4}
\put(17.3,7.6){\small 3} \put(17.3,5){\small 2}
\put(18,3.3){\small 1}

\put(13,0.5){\mbox{${\cal S}_3=\Sigma_{xy}S^{\,1}_5$}}


\put(22,3){\line(3,5){6}} \put(22,3){\line(1,1){4}}
\put(22,3){\line(3,1){6}} \put(22,3){\line(1,0){12}}

\put(34,3){\line(-3,5){6}} \put(34,3){\line(-1,1){4}}
\put(34,3){\line(-3,1){6}} \put(28,5){\line(1,1){2}}
\put(28,5){\line(-1,1){2}} \put(26,7){\line(1,0){4}}
\put(26,7){\line(1,3){2}} \put(30,7){\line(-1,3){2}}

\put(26,12.5){$x_1$} \put(27.8,3.8){$y_1$} \put(20.7,3.7){$y_3$}
\put(29.5,5.5){$x_3$} \put(34.4,3){$y_2$} \put(25.5,5.5){$x_2$}

\put(23,0.5){\mbox{${\cal S}_4=S^{\,0}_2\ast S^{\,0}_2 \ast
S^{\,0}_2$}}


\put(39,3){\line(-3,5){6}} \put(39,3){\line(0,1){8}}
\put(39,3){\line(3,5){6}} \put(39,3){\line(-1,3){2}}
\put(39,3){\line(1,3){2}} \put(33,13){\line(1,0){12}}

\put(39,8){\line(-2,1){2}} \put(39,8){\line(2,1){2}}
\put(37,9){\line(-1,1){4}} \put(41,9){\line(1,1){4}}
\put(39,11){\line(-1,-1){2}} \put(39,11){\line(1,-1){2}}
\put(39,11){\line(-3,1){6}} \put(39,11){\line(3,1){6}}

\put(32,12.5){\small 1} \put(45.5,12.5){\small 2}
\put(40.9,9.6){\scriptsize 3} \put(39.5,8.7){\scriptsize 4}
\put(36.8,9.6){\scriptsize 5} \put(38.8,11.5){$x$}
\put(40,3.3){$y$}

\put(36.5,0.5){\mbox{${\cal S}_5=S^{\,0}_2 \ast S^{\,1}_5$}}

\end{picture}

\setlength{\unitlength}{2.5mm}

\begin{picture}(47,14)(0,0)

\thicklines

\put(7,3){\line(-3,5){6}} \put(7,3){\line(0,1){8}}
\put(7,3){\line(3,5){6}} \put(1,13){\line(1,0){12}}
\put(1,13){\line(3,-4){6}} \put(1,13){\line(1,-1){6}}
\put(1,13){\line(3,-2){6}} \put(1,13){\line(3,-1){6}}
\put(13,13){\line(-3,-4){6}} \put(13,13){\line(-1,-1){6}}
\put(13,13){\line(-3,-2){6}} \put(13,13){\line(-3,-1){6}}

\put(0.2,12){$x$} \put(13.2,12){$y$} \put(7.7,3){\small 1}

\put(7.3,4.6){\scriptsize 2} \put(7.3,6.4){\scriptsize  3}
\put(7.3,8.3){\scriptsize  4} \put(7,11.5){\scriptsize 5}

\put(3,0.5){\mbox{${\cal S}_6 =\Sigma_{xy}S^{\,1}_6$}}


\thicklines

\put(12,3){\line(3,5){6}} \put(12,3){\line(1,1){6}}
\put(12,3){\line(2,1){4}} \put(12,3){\line(3,1){6}}
\put(12,3){\line(1,0){12}} \put(16,5){\line(1,0){4}}
\put(16,5){\line(1,2){2}} \put(20,5){\line(-1,2){2}}
\put(24,3){\line(-3,5){6}} \put(24,3){\line(-1,1){6}}
\put(24,3){\line(-2,1){4}} \put(24,3){\line(-3,1){6}}
\put(18,5){\line(0,1){8}}

\put(18.5,9.2){$x$} \put(11,3){\small 1} \put(16.5,12){\small 2}
\put(24.4,3){\small 3}  \put(15.3,5.2){\scriptsize  6}
 \put(20.2,5.2){\scriptsize 4} \put(17.8,3.7){\scriptsize 5}

\put(18,0.5){\mbox{${\cal S}_7$}}


\put(29,3){\line(-3,5){6}} \put(29,3){\line(3,5){6}}
\put(29,3){\line(-1,3){2}} \put(29,3){\line(1,3){2}}
\put(29,3){\line(-1,5){1}} \put(29,3){\line(1,5){1}}
\put(23,13){\line(1,0){12}} \put(27,9){\line(1,0){4}}
\put(28,8){\line(1,0){2}} \put(28,8){\line(-1,1){5}}
\put(30,8){\line(1,1){5}} \put(27,9){\line(2,1){8}}
\put(28,8){\line(3,1){3}}

\put(27.4,3){$x$} \put(22.3,12){\small 1} \put(35.2,12){\small 6}
\put(27,9.5){\scriptsize 2} \put(28.5,6.9){\scriptsize 3}
\put(30.7,7){\scriptsize 4} \put(30.6,9.3){\scriptsize 5}

\put(28.5,0.5){\mbox{${\cal S}_8$}}


\put(34,3){\line(3,5){6}} \put(34,3){\line(1,1){4}}
\put(34,3){\line(3,1){6}} \put(34,3){\line(1,0){12}}

\put(46,3){\line(-3,5){6}} \put(46,3){\line(-1,1){4}}
\put(46,3){\line(-3,1){6}} \put(40,5){\line(1,1){2}}
\put(40,5){\line(-1,1){2}} \put(38,7){\line(1,0){4}}
\put(38,7){\line(1,3){2}} \put(42,7){\line(-1,3){2}}

\put(38,7){\line(1,1){2}} \put(42,7){\line(-1,1){2}}
\put(40,9){\line(0,1){4}}

\put(33,3){\small 6} \put(46.4,3){\small 5}
\put(39.8,5.6){\scriptsize 4} \put(39.6,7.6){$x$}
\put(37.3,7.2){\scriptsize 2} \put(42.3,7.2){\scriptsize 3}
\put(38.5,12){\small 1}

\put(39.5,0.5){\mbox{${\cal S}_9$}}

\end{picture}

The following result (which we need later) follows from the
classification of combinatorial 2-spheres on $\leq 7$ vertices
(e.g., see \cite{a1, bd2}).

\begin{prop}$\!\!${\bf .} \label{s27}
Let $K$ be an $n$-vertex combinatorial $2$-sphere. If $n \leq 7$
then $K$ is isomorphic to ${\cal S}_1 \dots, {\cal S}_8$ or ${\cal
S}_9$ above.
\end{prop}

If $\kappa_{\alpha}$ is a proper bistellar $i$-move on a pure
simplicial complex $M$ and ${\rm lk}_M(\alpha) =
S^{i-1}_{i+1}(\beta)$ then $\beta$ is a removable $i$-face of
$\kappa_{\alpha}(M)$ (with ${\rm lk}_{\kappa_{\alpha}(M)}(\beta)
= S^{d-i - 1}_{d - i + 1}(\alpha)$) and $\kappa_{\beta} :
\kappa_{\alpha}(M) \mapsto M$ is a bistellar $(d - i)$-move. For
a vertex $u$, if ${\rm lk}_M(u) = S^{d-1}_{d+1}(\beta)$ then the
bistellar $d$-move $\kappa_{\{u\}} : M \mapsto \kappa_{\{u\}}(M)
= N$ deletes the vertex $u$ (we also say that $N$ is obtained
from $M$ by {\em collapsing} the vertex $u$). The operation
$\kappa_{\beta} : N \mapsto M$ is called a bistellar {\em
$0$-move}. We also say that $M$ is obtained from $N$ by {\em
starring} the vertex $u$ in the facet $\beta$ of $N$.

\setlength{\unitlength}{3mm}

\begin{picture}(45,12)(0,-2.5)


\thicklines

\put(2,5){\line(4,3){4}} \put(2,5){\line(3,-1){3}}
\put(2,5){\line(1,-1){4}} \put(8,5){\line(-2,3){2}}
\put(8,5){\line(-3,-1){3}} \put(8,5){\line(-1,-2){2}}
\put(5,4){\line(1,4){1}} \put(5,4){\line(1,-3){1}}

\put(14,5){\line(4,3){4}} \put(14,5){\line(3,-1){3}}
\put(14,5){\line(1,-1){4}} \put(20,5){\line(-2,3){2}}
\put(20,5){\line(-3,-1){3}} \put(20,5){\line(-1,-2){2}}
\put(17,4){\line(1,4){1}} \put(17,4){\line(1,-3){1}}
\put(18,1){\line(0,1){3.1}} \put(18,4.6){\line(0,1){3.4}}

\put(9,6){\vector(1,0){4}} \put(13,4){\vector(-1,0){4}}
\put(32,6){\vector(1,0){4}} \put(36,4){\vector(-1,0){4}}

\put(24,4){\line(1,1){4}} \put(24,4){\line(2,-1){4}}
\put(28,2){\line(0,1){6}} \put(31,4){\line(-3,-2){3}}
\put(31,4){\line(-3,4){3}}

\put(37,4){\line(1,1){4}} \put(37,4){\line(2,-1){4}}
\put(41,2){\line(0,1){6}} \put(44,4){\line(-3,-2){3}}
\put(44,4){\line(-3,4){3}}

\thinlines

\put(2,5){\line(1,0){2.9}} \put(8,5){\line(-1,0){2.4}}

\put(14,5){\line(1,0){2.9}} \put(20,5){\line(-1,0){1.7}}

\put(17.5,4.7){\mbox{-}}

\put(24,4){\line(1,0){3.6}} \put(28.4,4){\line(1,0){2.6}}

\put(37,4){\line(1,0){3}} \put(41.4,4){\line(1,0){2.6}}

\put(40,5){\line(1,-3){1}} \put(40,5){\line(-3,-1){3}}
\put(40,5){\line(1,3){1}} \put(44,4){\line(-4,1){2.6}}

\put(40.2,4.8){\mbox{$.$}} \put(40.65,4.752){\mbox{$.$}}
\put(40.5,3.948){\mbox{$.$}}

\put(9,7){\mbox{{\small $1$-move}}} \put(9.3,2.5){\mbox{{\small
$2$-move}}}

\put(32,7){\mbox{{\small $0$-move}}} \put(32.3,2.5){\mbox{{\small
$3$-move}}}

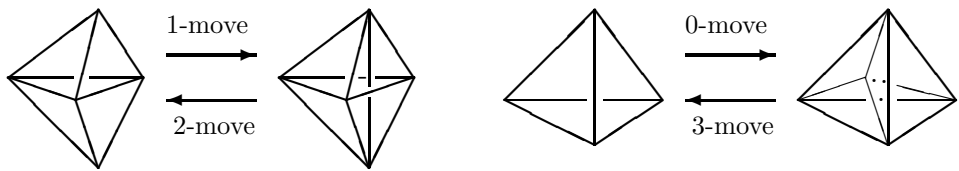
\put(14,-1.5){\mbox{Bistellar moves in dimension $3$}}

\end{picture}

If $M$ is a 3-pseudomanifold and $\kappa_{\alpha} : M \mapsto N$
is a bistellar 1-move then, from the definition, $(f_0(N)$,
$f_1(N), f_2(N), f_3(N)) = (f_0(M), f_1(M) + 1, f_2(M) + 2,
f_3(M) + 1)$ and $\deg_{N}(v) \geq \deg_M(v)$ for any vertex $v$.

Two simplicial complexes $K$ and $L$ are called {\em bistellar
equivalent} (denoted by $K\sim L$) if there exists a finite
sequence of bistellar moves leading from $K$ to $L$. Let
$\kappa_{\alpha}$ be a proper bistellar $i$-move and ${\rm
lk}_M(\alpha) = S^{i-1}_{i+1}(\beta)$. If $K_1$ is obtained from
$K$ by starring (\cite{bd2}) a new vertex in $\alpha$ and $K_2$
is obtained from $\kappa_{\alpha}(K)$ by starring a new vertex in
$\beta$ then $K_1$ and $K_2$ are isomorphic. Thus, if $K \sim L$
then $K \approx L$. Conversely, it was shown in \cite{p}, that if
two combinatorial manifolds are combinatorially equivalent then
they are bistellar equivalent.

Let $\tau\subset\sigma$ be two faces of a simplicial complex $K$.
We say that $\tau$ is a {\em free face} of $\sigma$ if $\sigma$ is
the only face of $K$ which properly contains $\tau$. (It follows
that $\dim(\sigma)-\dim(\tau)=1$ and $\sigma$ is a maximal simplex
in $K$.) If $\tau$ is a free face of $\sigma$ then $K^{\,\prime}
:= K \setminus \{\tau, \sigma\}$ is a simplicial complex. We say
that there is an {\em elementary collapse} of $K$ to
$K^{\,\prime}$. We say $K$ {\em collapses} to $L$ and write
$K\scoll L$ if there exists a sequence $K=K_0$, $K_1, \dots$,
$K_n=L$ of simplicial complexes such that there is an elementary
collapse of $K_{i-1}$ to $K_{i}$ for $1\leq i\leq n$. If $L$
consists of a 0-simplex (a vertex) we say that $K$ is {\em
collapsible} and write $K\scoll 0$.

Suppose $P^{\,\prime}\subseteq P$ are polyhedra and $P =
P^{\,\prime} \cup B$, where $B$ is a   pl $k$-ball (for some
$k\geq 1$). If $P^{\,\prime} \cap B$ is a   pl $(k-1)$-ball then
we say that there is an {\em elementary collapse} of $P$ to
$P^{\,\prime}$. We say that $P$ collapses to $Q$ and write $P\coll
Q$ if there exists a sequence $P = P_0, P_1, \dots, P_n = Q$ of
polyhedra such that there is an elementary collapse of $P_{i-1}$
to $P_{i}$ for $1\leq i\leq n$. For two simplicial complexes $K$
and $L$, if $K\scoll L$  then clearly $|K| \coll |L|$. The
following is a consequence of the Simplicial Neighbourhood
Theorem (see \cite{bd5}).

\begin{prop}$\!\!${\bf .} \label{snt}
Let $\sigma$ be a facet of a connected combinatorial $d$-manifold
$X$.  Put $L = C(\sigma, X)$, the simplicial complement of
$\sigma$ in $X$. Also, let $Y = X \setminus \{\sigma\}$. Then
\begin{enumerate}
     \item[$(a)$] $|Y| \coll |L|$.
     \item[$(b)$] If, further, $L$ is collapsible then $X$ is a
     combinatorial sphere.
     \end{enumerate}
\end{prop}

\section{Proof of Theorem 1}

For $n\geq 4$, by an $S^{\,2}_n$ we mean a combinatorial 2-sphere
on $n$ vertices. If a combinatorial 3-manifold has at most 9
vertices then it is connected and hence is a 3-pseudomanifold.

\begin{lemma}$\!\!${\bf .} \label{l3.1}
Let $N$ be an $n$-vertex combinatorial $3$-manifold with minimum
vertex-degree $k \leq n-2$ and $n \leq 9$. Let $u$ be a vertex of
degree $k$ in $N$. Then there exists a bistellar $1$-move
$\kappa_{\beta} : N \mapsto \widetilde{N}$ such that
$\deg_{\widetilde{N}}(u) = k+1$.
\end{lemma}

\noindent {\bf Proof.} Let $X = {\rm lk}_N(u) = S^{\,2}_k$. If $k
= 4$ then $X = S^{\,2}_4(\{a, b, c, d\})$ for some $a, b, c, d\in
V(N)$. Let $\beta = abc$. Suppose ${\rm lk}_N(\beta) = \{u, x\}$.
If $x = d$ then the induced subcomplex $K = N[\{u, a, b, c, d\}]$
is a 3-pseudomanifold. Since $n \geq 6$, $K$ is a proper
subcomplex of $N$. This is not possible. Thus, $x \neq d$ and
hence $ux$ is a non-edge in $N$. So, $\kappa_{\beta}$ is a
bistellar 1-move, as required. So, let $k\geq 5$.

Suppose the result is false. Let ${\cal B}$ denote the collection
of all facets $B\in N$ such that $u\not\in B$ and $B$ contains a
triangle of $X$. Then ${\cal B}$ is a set of 4-sets satisfying
$(a)$ each element of ${\cal B}$ is contained in $V(X)$, $(b)$
each triangle of $X$ is contained in a unique member of ${\cal
B}$, and $(c)$ each member of ${\cal B}$ contains one or two
triangles of $X$. (Indeed, if $B \in {\cal B}$ is not contained in
$V(X)$ and $\beta \subseteq B$ is a triangle of $X$, then
$\kappa_{\beta}$ is a 1-move on $N$ which increases $\deg(u)$,
contrary to our assumption that $N$ does not admit such a move.
This proves (a). Now, (b) is immediate since $N$ is a
pseudomanifold. If $B \in {\cal B}$ contains three triangles of
$X$ then these three triangles have a common vertex $x$, and
$\deg_N(x) = 4$, contrary to our assumption that the minimum
vertex-degree of $N$ is $\geq 5$. This proves (c)).

Let $G$ denote the graph whose vertices are the 4-subsets of
$V(X)$ containing 1 or 2 triangles of $X$. Two vertices are
adjacent if the corresponding 4-subsets have a triangle of $X$ in
common. It follows that ${\cal B}$ is a maximal coclique in $G$.
So, we look for the maximal cocliques of $G$ for each admissible
choice of $X$.

In case $k = 5$, $X$ is of the form $S^{\,0}_2(xy) \ast
S^{\,1}_3(abc)$. In this case, $G$ has a unique maximal coclique
$ C = \{xyab, xyac, xybc\}$. If ${\cal B} = C$, then $N$ contains
a proper 3-dimensional sub-pseudomanifold $S^{\,1}_3(uxy) \ast
S^{\,1}_3(abc)$, a contradiction. So, $k\geq 6$.

Let $k =6$. Then $X$ is isomorphic to ${\cal S}_3$ or ${\cal S}_4$
(defined in Example \ref{e2}). Consider the case $X = {\cal
S}_3$. Then the vertices of $G$ are $v_{1} = x123$, $v_{2} =
y123$, $v_{3} = x234$, $v_{4} = y234$, $v_{5} = xy14$, $v_{6} =
xy23$, $v_{7} = x124$, $v_{8} = y124$, $v_{9} = x134$, $v_{10} =
y134$, $v_{11} = xy13$, $v_{12} = xy24$ and

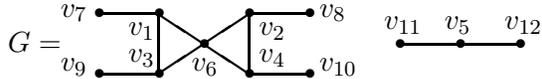
\begin{picture}(40,5)(0,0)
\setlength{\unitlength}{2mm}

\put(5.7,1){$_{\bullet}$} \put(5.7,5){$_{\bullet}$}
\put(9.7,1){$_{\bullet}$} \put(9.7,5){$_{\bullet}$}
\put(12.7,3){$_{\bullet}$} \put(15.7,1){$_{\bullet}$}
\put(15.7,5){$_{\bullet}$} \put(19.7,1){$_{\bullet}$}
\put(19.7,5){$_{\bullet}$} \put(25.7,3){$_{\bullet}$}
\put(29.7,3){$_{\bullet}$} \put(33.7,3){$_{\bullet}$}

\thicklines

\put(6,1){\line(1,0){4}} \put(6,5){\line(1,0){4}}
\put(10,5){\line(3,-2){3}} \put(10,1){\line(3,2){3}}
\put(13,3){\line(3,-2){3}} \put(13,3){\line(3,2){3}}
\put(16,1){\line(1,0){4}} \put(16,5){\line(1,0){4}}
\put(10,1){\line(0,1){4}} \put(16,1){\line(0,1){4}}
\put(26,3){\line(1,0){8}}

\put(0,2.5){$G=$}

\put(3.5,1.2){$v_{9}$} \put(3.5,4.8){$v_{7}$} \put(8,1.8){$v_{3}$}
\put(8,3.6){$v_{1}$} \put(12.2,1.2){$v_{6}$}
\put(16.7,1.8){$v_{4}$} \put(20.7,1.2){$v_{10}$}
\put(16.7,3.6){$v_{2}$} \put(20.7,4.8){$v_{8}$}
\put(25,4){$v_{11}$} \put(29,4){$v_{5}$} \put(33,4){$v_{12}$}

\end{picture}

Note that ${\rm Aut}(X) \cong \ZZ_2 \times \ZZ_2$, generated by
$(x, y)$ and $(1, 4)(2, 3)$. Up to this automorphism group, the
maximal cocliques of $G$ are $C_1 = \{v_{6}, v_{7}, v_{8}, v_{9},
v_{10}, v_{11}, v_{12}\}$, $C_2 = \{v_{5}, v_{6}, v_{7}, v_{8},
v_{9}, v_{10}\}$, $C_3 = \{v_{3}, v_{4}, v_{7}, v_{8}, v_{11},
v_{12}\}$, $C_4 = \{v_{2}, v_{3}, v_{7}, v_{10}, v_{11},
v_{12}\}$, $C_5 = \{v_{3}, v_{4}, v_{5}, v_{7}, v_{8}\}$ and $C_6
= \{v_{2}, v_{3}, v_{5}, v_{7}, v_{10}\}$.

If ${\cal B} = C_i$, for $1\leq i\leq 4$, then the available
portion of ${\rm lk}_N(x)$ can not be completed to a 2-sphere, a
contradiction. If ${\cal B} = C_5$ then ${\rm lk}_N(1) =
S^{\,0}_2(u4) \ast S^{\,1}_3(xyz)$, so that $\deg_N(1) = 5 < k$,
a contradiction. If ${\cal B} = C_6$ then ${\rm lk}(x)$ and ${\rm
lk}(y)$ are $S^{\,2}_6$'s with vertex-sets $\{u, 1, 2, 3, 4, y\}$
and $\{u, 1, 2, 3, 4, x\}$ respectively. It follows that the
remaining (one or two) vertices can only be joined to each other
and with $1, 2, 3, 4$. Then these vertices have degree $\leq 5$,
a contradiction.

Next assume that $X= {\cal S}_4$. Up to automorphism of $X$, there
are two maximal cocliques of $G$, namely, $C_1 = \{x_1y_1x_2x_3,
x_1y_1x_2y_3, x_1y_1y_2x_3, x_1y_1y_2y_3\}$ and $C_2 =
\{x_1y_1x_2x_3, x_1y_1x_2y_3$, $x_1y_2x_3y_3, y_1y_2x_3y_3\}$. If
${\cal B} = C_1$ or $C_2$, then ${\rm lk}_N(x_2) = S^{\,0}_2(x_1
y_1) \ast S^{\,1}_3(ux_3y_3)$ and hence $\deg_N(x_2) = 5 <k$, a
contradiction.

Thus $k =7$ and $n=9$. Let $uv$ be a non-edge. Then $f_1(N) \leq
35$ and hence $f_3(N) \leq 26$. Since there are 10 facets through
$u$ and 10 facets through $v$, it follows that $\#({\cal B}) \leq
6$. Since there are 10 triangles in $X$ and each member of ${\cal
B}$ contains at most two triangles of $X$, $\#({\cal B}) \geq 5$.
Thus ${\cal B}$ is a maximal coclique of $G$ of size 5 or 6.

Since $X$ has 7 vertices, $X$ is isomorphic to ${\cal S}_5, \dots,
{\cal S}_8$ or ${\cal S}_9$ (of Example \ref{e2}).

Consider the case, $X= {\cal S}_5$. Then the vertices of $G$ are
$v_{1} = x124$, $v_{2} = x235$, $v_{3} = x134$, $v_{4} = x245$,
$v_{5} = x135$, $v_{6} = y124$, $v_{7} = y235$, $v_{8} = y134$,
$v_{9} = y245$, $v_{10} = y135$, $v_{11} = xy12$, $v_{12} =
xy23$, $v_{13} = xy34$, $v_{14} = xy45$, $v_{15} = xy15$, $v_{16}
= x123$, $v_{17} = x234$, $v_{18} = x345$, $v_{19} = x145$,
$v_{20} = x125$, $v_{21} = y123$, $v_{22} = y234$, $v_{23} =
y345$, $v_{24} = y145$ and $v_{25} = y125$. Note that ${\rm
Aut}(X) \cong \ZZ_2 \times D_{10}$, generated by $(x, y)$, $(1, 2,
3, 4, 5)$ and $(1, 2)(3, 5)$. Up to this automorphism group,
there are only 11 maximal 6-cocliques and 3 maximal 5-cocliques.
These are $C_1 = \{v_{1}, v_{2}, v_{13}, v_{14}, v_{15},
v_{21}\}$, $C_2 = \{v_{1}, v_{2}, v_{13}, v_{19}, v_{21},
v_{24}\}$, $C_3 = \{v_{1}, v_{2}, v_{15}, v_{18}, v_{21},
v_{23}\}$, $C_4 = \{v_{1}, v_{3}, v_{12}, v_{19}, v_{23},
v_{25}\}$, $C_5 = \{v_{1}, v_{6}, v_{12}, v_{13}, v_{14},
v_{15}\}$, $C_6 = \{v_{1}, v_{6}, v_{12}, v_{13}, v_{19},
v_{24}\}$, $C_7 = \{v_{1}, v_{6}, v_{12}, v_{15}, v_{18},
v_{23}\}$, $C_8 = \{v_{1}, v_{6}, v_{17}, v_{19}, v_{22},
v_{24}\}$, $C_9 = \{v_{1}, v_{7}, v_{17}, v_{19}, v_{23},
v_{25}\}$, $C_{10} = \{v_{1}, v_{8}, v_{14}, v_{15}, v_{17},
v_{21}\}$, $C_{11} = \{v_{1}, v_{8}, v_{17}, v_{19}, v_{21},
v_{24}\}$, $C_{12} = \{v_{11}, v_{12}, v_{13}, v_{14}, v_{15}\}$,
$C_{13} = \{v_{11}, v_{12}, v_{13}, v_{19}, v_{24}\}$ and $C_{14}
= \{v_{11}, v_{17}, v_{19}, v_{22}, v_{24}\}$. If ${\cal B} =
C_i$, for $i = 1, \dots, 5$ or 7, then $x$ becomes a singular
vertex, a contradiction. If ${\cal B} = C_9$, then $x$ becomes a
vertex of degree 6, a contradiction. For $i\in \{6, 8, 10, 11,
12, 13, 14\}$, if ${\cal B} = C_{i}$ then $5$ becomes a vertex of
degree 5, a contradiction.

Now, let $X= {\cal S}_6$. Then the vertices of $G$ are $v_{1} =
xy13$, $v_{2} = xy14$, $v_{3} = xy25$, $v_{4} = xy35$, $v_{5} =
x124$, $v_{6} = y124$, $v_{7} = x125$, $v_{8} = y125$, $v_{9} =
x134$, $v_{10} = y134$, $v_{11} = x235$, $v_{12} = y235$, $v_{13}
= x245$, $v_{14} = y245$, $v_{15} = x145$, $v_{16} = y145$,
$v_{17} = xy15$, $v_{18} = xy23$, $v_{19} = xy34$, $v_{20} =
x123$, $v_{21} = y123$, $v_{22} = x234$, $v_{23} = y234$, $v_{24}
= x345$, $v_{25} = y345$. Note that ${\rm Aut}(X) \cong \ZZ_2
\times \ZZ_2$, generated by $(x, y)$ and $(1, 5)(2, 4)$. Up to
${\rm Aut}(X)$, there are 11 maximal 6-cocliques and one maximal
5-coclique. These are $C_1 = \{v_{1}, v_{3}, v_{20}, v_{21},
v_{24}, v_{25}\}$, $C_2 = \{v_{1}, v_{4}, v_{20}, v_{21}, v_{24},
v_{25}\}$, $C_3 = \{v_{2}, v_{3}, v_{20}, v_{21}, v_{24},
v_{25}\}$, $C_4 = \{v_{5}, v_{6}, v_{17}, v_{18}, v_{24},
v_{25}\}$, $C_5 = \{v_{5}, v_{8}, v_{17}, v_{18}, v_{24},
v_{25}\}$, $C_6 = \{v_{5}, v_{11}, v_{17}, v_{21}, v_{24},
v_{25}\}$, $C_7 = \{v_{5}, v_{13}, v_{17}, v_{21}, v_{22},
v_{25}\}$, $C_8 = \{v_{5}, v_{15}, v_{17}, v_{21}, v_{22},
v_{25}\}$, $C_9 = \{v_{7}, v_{8}, v_{17}, v_{18}, v_{24},
v_{25}\}$, $C_{10} = \{v_{7}, v_{11}, v_{17}, v_{21}, v_{24}$,
$v_{25}\}$, $C_{11} = \{v_{7}, v_{15}, v_{17}, v_{21}, v_{22},
v_{25}\}$, $C_{12} = \{v_{17}, v_{20}, v_{21}, v_{24}, v_{25}\}$.
For $i\in \{1, 3, 4, 5, 6, 7, 9$, $11\}$, if ${\cal B} = C_i$ then
$x$ becomes a singular vertex, a contradiction. If ${\cal B} =
C_2$ or $C_{12}$ then $\deg(2) = 5$, a contradiction. If ${\cal
B} = C_8$ or $C_{10}$ then ${\rm lk}(x)$ is an $S^{\,2}_7$ with
vertex-set $\{u, y, 1, \dots, 5\}$. So, $uv$ and $xv$ are
non-edges and hence $\deg(v) \leq 6$, a contradiction.

Consider the case, $X= {\cal S}_7$. Then the vertices of $G$ are
$v_{1} = x135$, $v_{2} = x124$, $v_{3} = x125$, $v_{4} = x235$,
$v_{5} = x236$, $v_{6} = x346$, $v_{7} = x134$, $v_{8} = x145$,
$v_{9} = x245$, $v_{10} = x256$, $v_{11} = x356$, $v_{12} = x136$,
$v_{13} = x146$, $v_{14} = 1234$, $v_{15} = 1236$, $v_{16} =
2345$, $v_{17} = 3456$, $v_{18} = 1256$, $v_{19} = 1456$, $v_{20}
= x126$, $v_{21} = x234$, $v_{22} = x456$, $v_{23} = 1235$,
$v_{24} = 1345$, $v_{25} = 1356$. In this case, ${\rm Aut}(X)
\cong D_6$, generated by $(1, 3, 5)(2, 4, 6)$ and $(1, 3)(4, 6)$.
There is no maximal 5-coclique and up to ${\rm Aut}(X)$, there are
only 3 maximal 6-cocliques. These are $C_1 = \{v_{16}, v_{18},
v_{20}, v_{21}, v_{22}, v_{23}\}$, $C_2 = \{v_{16}, v_{19},
v_{20}, v_{21}, v_{22}, v_{23}\}$, $C_3 = \{v_{17}, v_{19},
v_{20}, v_{21}, v_{22}, v_{23}\}$. If ${\cal B} = C_1$ or $C_2$
then 3 becomes a vertex of degree 6, a contradiction. If ${\cal B}
= C_3$ then 5 becomes a singular vertex, a contradiction.

Consider the case, $X= {\cal S}_8$. Then the vertices of $G$ are
$v_{1} = x124$, $v_{2} = x125$, $v_{3} = x236$, $v_{4} = x134$,
$v_{5} = x346$, $v_{6} = x145$, $v_{7} = x245$, $v_{8} = x356$,
$v_{9} = x136$, $v_{10} = x146$, $v_{11} = 1236$, $v_{12} = 1246$,
$v_{13} = 2456$, $v_{14} = 1235$, $v_{15} = 1345$, $v_{16} =
3456$, $v_{17} = x123$, $v_{18} = x234$, $v_{19} = x456$, $v_{20}
= x156$, $v_{21} = x235$, $v_{22} = x256$, $v_{23} = 1256$,
$v_{24} = 2345$, $v_{25} = 2356$. Here, ${\rm Aut}(X) \cong
\ZZ_2$, generated by $(1, 4)(2, 5)(3, 6)$. There is no maximal
5-coclique and up to ${\rm Aut}(X)$, there are 10 maximal
6-cocliques. These are $C_1 = \{v_{4}, v_{10}, v_{17}, v_{19},
v_{23}, v_{24}\}$, $C_2 = \{v_{5}, v_{9}, v_{17}, v_{19}, v_{23},
v_{24}\}$, $C_3 = \{v_{4}, v_{9}, v_{17}, v_{19}, v_{23},
v_{24}\}$, $C_4 = \{v_{1}, v_{6}, v_{18}, v_{20}, v_{23},
v_{24}\}$, $C_5 = \{v_{2}, v_{7}, v_{18}, v_{20}, v_{23},
v_{24}\}$, $C_6 = \{v_{1}, v_{7}, v_{18}, v_{20}, v_{23},
v_{24}\}$, $C_7 = \{v_{4}, v_{6}, v_{17}, v_{20}, v_{23},
v_{24}\}$, $C_8 = \{v_{4}, v_{7}, v_{17}, v_{20}, v_{23},
v_{24}\}$, $C_9 = \{v_{5}, v_{6}, v_{17}, v_{20}, v_{23},
v_{24}\}$, $C_{10} = \{v_{5}, v_{7}, v_{17}, v_{20}, v_{23},
v_{24}\}$. If ${\cal B} = C_i$, for $i = 1, 2, 4, 5$ or 7, then
${\rm lk}(x)$ is an $S^{\,2}_7$ with vertex-set $\{u, y, 1, \dots,
5\}$. So, $uv$ and $xv$ are non-edges and hence $\deg(v) \leq 6$,
a contradiction. If ${\cal B} = C_i$, for $i = 3, 6, 8, 9$ or 10
then $x$ becomes a singular vertex, a contradiction.

Consider the case, $X= {\cal S}_9$. Then the vertices of $G$ are
$v_{1} = x124$, $v_{2} = x125$, $v_{3} = x235$, $v_{4} = x236$,
$v_{5} = x134$, $v_{6} = x136$, $v_{7} = x156$, $v_{8} = x246$,
$v_{9} = x345$, $v_{10} = x456$, $v_{11} = 1234$, $v_{12} = 1235$,
$v_{13} = 1236$, $v_{14} = x126$, $v_{15} = x234$, $v_{16} =
x135$, $v_{17} = 1246$, $v_{18} = 1256$, $v_{19} = 2345$, $v_{20}
= 2346$, $v_{21} = 1345$, $v_{22} = 1356$, $v_{23} = 1456$,
$v_{24} = 2456$, $v_{25} = 3456$. Here, ${\rm Aut}(X) \cong D_6$,
generated by $(1, 2, 3)(4, 5, 6)$ and $(1, 2)(4, 5)$. There is no
maximal 5-coclique and up to ${\rm Aut}(X)$, there are 13 maximal
6-cocliques. These are $C_1 = \{v_{1}, v_{5}, v_{15}, v_{17},
v_{22}, v_{25}\}$, $C_2 = \{v_{1}, v_{6}, v_{15}, v_{17}, v_{22},
v_{25}\}$, $C_3 = \{v_{2}, v_{5}, v_{15}, v_{17}, v_{22},
v_{25}\}$, $C_4 = \{v_{2}, v_{6}, v_{15}, v_{17}, v_{22},
v_{25}\}$, $C_5 = \{v_{1}, v_{5}, v_{15}, v_{17}, v_{21},
v_{23}\}$, $C_6 = \{v_{1}, v_{6}, v_{15}, v_{17}, v_{21},
v_{23}\}$, $C_7 = \{v_{2}, v_{5}, v_{15}, v_{17}, v_{21},
v_{23}\}$, $C_8 = \{v_{2}, v_{6}, v_{15}, v_{17}, v_{21},
v_{23}\}$, $C_9 = \{v_{8}, v_{9}, v_{14}, v_{15}, v_{16},
v_{23}\}$, $C_{10} = \{v_{1}, v_{9}, v_{15}, v_{16}, v_{17},
v_{23}\}$, $C_{11} = \{v_{2}, v_{9}, v_{15}, v_{16}, v_{17},
v_{23}\}$, $C_{12} = \{v_{1}, v_{7}, v_{15}, v_{16}, v_{17},
v_{25}\}$, $C_{13} = \{v_{2}, v_{7}, v_{15}, v_{16}, v_{17},
v_{25}\}$. If ${\cal B} = C_3, C_4, C_8$ or $C_{13}$, then $5$
becomes a singular vertex, a contradiction. If ${\cal B} = C_1$ or
$C_5$ then $\deg(x) = 5$, a contradiction. If ${\cal B} = C_{2},
C_6, C_9$ or $C_{12}$,  then $\deg(2) = 6$, a contradiction. If
${\cal B} = C_{7}, C_{10}$ or $C_{11}$, then $\deg(3) = 6$, a
contradiction. This completes the proof. \hfill $\Box$

\bigskip

\noindent {\bf Proof of Theorem 1.} Let $M$ be a $9$-vertex
combinatorial 3-manifold. Then, by Lemma \ref{l3.1}, there exists
a sequence of bistellar 1-moves $\kappa_{\beta_1}, \dots,
\kappa_{\beta_l}$ such that $N := \kappa_{\beta_l}(\cdots
(\kappa_{\beta_1}(M)))$ is neighbourly. Since $f_1(M) \geq
4\times 9 - 10 = 26$ (see \cite{w}) and each bistellar 1-move
produces an edge, $l \leq 10$. Let $\alpha_i = {\rm
lk}_{\kappa_{\beta_{i - 1}}(\cdots
(\kappa_{\beta_1}(M)))}(\beta_i)$ for $2\leq i\leq l$ and
$\alpha_1 = {\rm lk}_{M}(\beta_1)$. Then $M =
\kappa_{\alpha_1}(\cdots (\kappa_{\alpha_l}(N)))$. This completes
the proof. \hfill  $\Box$

\begin{remark}$\!\!${\bf .}
{\rm Let $C^{\,3}_7$ be the cyclic $3$-sphere whose facets are
those $4$-subsets of the vertices of the 7-cycle $S^{\,1}_7(1
\cdots 7)$ on which the 7-cycle induces a subgraph with
even-sized components. Let $K$ be the simplicial complex
$C^{\,3}_7\setminus\{1234\}$. Let $K^{\prime}$ be the simplicial
complex on the vertex-set $\{1^{\prime}, \dots, 7^{\prime}\}$
isomorphic to $K$ (by the map $i\mapsto i^{\prime}$). Consider the
simplicial complex $M$ which is obtained from $K \sqcup
K^{\prime}$ by identifying $i$ with $i^{\prime}$ for $1\leq i\leq
4$. Then $M$ is a combinatorial $3$-sphere (connected sum of two
copies of $C^{\,3}_7$). The minimum vertex-degree in $M$ is 6. The
degree of the vertex $6$ in $M$ is $6$ with non-edges
$65^{\prime}$, $66^{\prime}$ and $67^{\prime}$. But, there is no
bistellar 1-move $\kappa_{\beta} : M \mapsto \kappa_{\beta}(M)$
such that $\deg_{\kappa_{\beta}(M)}(6) = 7$. So, Lemma \ref{l3.1}
is not true for $n=10$.}
\end{remark}

\begin{remark}$\!\!${\bf .}
{\rm Let $X$ be an $S^{\,2}_k$. Let $\alpha = \alpha(X)$ denote
the number of 4-subsets of $V(X)$ which contain one or two
triangles of $X$. While proving Lemma \ref{l3.1}, we noticed that
when $k = 6$, $\alpha(X) = 12$, and when $k = 7$, $\alpha(X) =
25$ : independent of the choice of $X$! This is no accident.
Indeed, we have $\alpha(S^{\,2}_k) = (k-2)(2k-9)$, regardless of
the choice of $S^{\,2}_k$. This may be proved by noting that any
two $S^{\,2}_k$ are equivalent by a sequence of proper bistellar
moves, and $\alpha$ is invariant under such moves. The explicit
value of $\alpha$ may then be computed by making a judicious
choice of $S^{\,2}_k$. }
\end{remark}

\section{Proof of Theorem 2}

{\sf Throughout this section $M^3_9$ will denote a fixed but
arbitrary neighbourly non-sphere combinatorial 3-manifold}. Thus,
$\chi(M^3_9) = 0$, and $f(M^3_9) = (9, 36, 54, 27)$.

\begin{lemma}$\!\!${\bf .} \label{l4.1}
The $f$-vector of the simplicial complement of any facet of
$M^3_9$ is either $(5, 10, 7, 1)$ or $(5, 10, 6, 0)$.
\end{lemma}

\noindent {\bf Proof.} Let $\sigma$ be a facet of $M^3_9$ and let
$(f_0, f_1, f_2, f_3)$ be the $f$-vector of the simplicial
complement $C(\sigma, M^3_9)$. By Proposition \ref{snt}, the
geometric carrier of the simplicial complex $M^3_9 \setminus
\{\sigma\}$ collapses to that of $C(\sigma, M^3_9)$. Since the
Euler characteristic is invariant under collapsing, we get
$\chi(C(\sigma, M^3_9)) = \chi(M^3_9 \setminus \{\sigma\}) = 1$.
Thus, $f_0 - f_1 + f_2 - f_3 =1$. Also, as $M^3_9$ is neighbourly,
$f_0 = 5$ and $f_1 = {5 \choose 2} = 10$. Hence $f_2 = f_3 +6$.
So, to complete the proof, it is sufficient to show that $f_3
\leq 1$.

Since $f_2 \leq {5 \choose 3} = 10$, it follows that $f_3\leq 4$.
Clearly, there are unique simplicial complexes with $f$-vectors
$(5, 10, 10, 4)$, $(5, 10, 9, 3)$ and $(5, 10, 8, 2)$, and all
these are collapsible. But, if $C(\sigma, M^3_9)$ was collapsible
then, by Proposition \ref{snt}, $M^3_9$ would be a sphere. So,
$f_3\leq 1$.  \hfill $\Box$

\begin{lemma}$\!\!${\bf .} \label{l4.2}
Let $\sigma_1$, $\sigma_2$ be two disjoint facets of $M^3_9$ and
let $x$ be the unique vertex of $M^3_9$ outside $\sigma_1\cup
\sigma_2$. Then the induced subcomplex of ${\rm lk}_{M^3_9}(x)$ on
$\sigma_1$ $($as well as on $\sigma_2)$ is an $S^{\,1}_3$
together with an isolated vertex.
\end{lemma}

\noindent {\bf Proof.} By Lemma \ref{l4.1}, the simplicial
complement $C(\sigma_2, M^3_9)$ of $\sigma_2$ has only one facet
(viz. $\sigma_1$) and seven triangles, four of which are the
triangles in $\sigma_1$. So, $C(\sigma_2, M^3_9)$ contains exactly
three triangles through $x$. Up to isomorphism, there are two
choices for these three triangles, one of which leads to a
collapsible complex $C(\sigma_2, M^3_9)$, which is not possible by
Proposition \ref{snt}. In the remaining case, we get the
situation as described in the lemma. \hfill $\Box$

\begin{lemma}$\!\!${\bf .} \label{l4.3}
Suppose each vertex of $M^3_9$ is in exactly two edges of degree
$3$. Then $M^3_9$ has an edge of degree $\geq 6$.
\end{lemma}

\noindent {\bf Proof.} Fix any facet $\sigma$ of $M^3_9$. Since
$M^3_9$ is neighbourly, the link of each vertex is an $S^{\,2}_8$
and hence has 12 triangles. Thus each vertex of $M^3_9$ is in 12
facets. Therefore, by the inclusion-exclusion principle, the
number of facets meeting $\sigma$ in at least one vertex is ${4
\choose 1}\times 12 - {\displaystyle \sum_{e\subset \sigma}
\deg(e)} + {4 \choose 3} \times 2 - {4 \choose 4} \times 1 = 55 -
{\displaystyle \sum_{e\subset \sigma} \deg(e)}$. Hence, by
subtraction, the number of facets disjoint from $\sigma$ is
${\displaystyle \sum_{e\subset \sigma} \deg(e)} - 28$. But, by
Lemma \ref{l4.1}, at most one facet can be disjoint from
$\sigma$. Hence
\begin{equation} \label{28/29}
\sum_{e\subset \sigma} \deg(e) = 29 \mbox{ or } 28
\end{equation}
according as there is a (necessarily unique) facet of $M^3_9$
disjoint from $\sigma$, or not. Here the sum is over all the six
edges of $M^3_9$ contained in the facet $\sigma$.

Now suppose, if possible, that all the edges of $M^3_9$ are of
degree $3, 4$ or 5. Then (\ref{28/29}) implies that any facet of
$M^3_9$ contains at most one edge of degree 3. Let $G$ denote the
graph with vertex-set $V(M^3_9)$ whose edges are precisely the
edges of degree 3 in $M^3_9$. By, our assumption, $G$ is a
9-vertex regular graph of degree 2, i.e., a disjoint union of
cycles. If $e=xy$ is any edge of $G$ then, putting $A = V({\rm
lk}(xy))$, we see that $A\cup \{x\}$ and $A\cup \{y\}$ are two
cocliques of size 4 in $G$. This is because no facet of $M^3_9$
contains more than one edge of $G$. But, we see by inspection
that the 9-cycle $S^{\,1}_9$ is the only 9-vertex union of cycles
in which there is such a pair of 4-cocliques corresponding to
every edge $e$. Thus, $G = S^{\,1}_9$. Also, for every edge
$e=xy$ of $S^{\,1}_9$, there is a unique set $A$ of vertices of
$S^{\,1}_9$ such that $A \cup \{x\}$ and $A \cup \{y\}$ are
4-cocliques of $S^{\,1}_9$.

This observation uniquely determines the link in $M^3_9$ of all
its degree 3 edges. Hence all the 27 distinct facets of $M^3_9$
are determined. But we now see that any two vertices at a
distance 2 in $S^{\,1}_9$ form an edge of degree 6 in $M^3_9$, a
contradiction. \hfill $\Box$

\begin{lemma}$\!\!${\bf .} \label{l4.4}
There is at least one pair of disjoint facets in $M^3_9$.
\end{lemma}

\noindent {\bf Proof.} Suppose not. Thus, any two facets of
$M^3_9$ intersect. Let $e$ be an edge of degree 7. Then $2\times
12 - 7 = 17$ facets intersect $e$ and hence $27 - 17 = 10$ facets
are disjoint from $e$. These facets are 4-sets in the heptagon
${\rm lk}(e)$, each of which meets all the edges of the heptagon.
But, one sees that the heptagon contains only seven such 4-sets,
a contradiction. So, $M^3_9$ has no edge of degree 7.

Next, let $e$ be an edge of degree 6. Let ${\rm lk}(e) =$
\setlength{\unitlength}{2mm}
\begin{picture}(5,2)(4,-0.5)

\thicklines

\put(3,0){\line(2,1){2}} \put(3,0){\line(2,-1){2}}
\put(5,1){\line(1,0){1.6}} \put(5,-1){\line(1,0){1.6}}
\put(6.6,1){\line(2,-1){2}} \put(6.6,-1){\line(2,1){2}}

\put(2.8,0.7){\scriptsize 1} \put(4.3,1.3){\scriptsize 2}
\put(6.5,1.3){\scriptsize 3} \put(8.2,0.7){\scriptsize 4}
\put(6.5,-0.6){\scriptsize 5} \put(4.7,-0.6){\scriptsize 6}
\end{picture} and let $x$ be the unique vertex outside $e
\cup \{1, \dots, 6\}$. Each of the $27 - 2\times 12 + 6 = 9$
facets disjoint from $e$ is a 4-set meeting all the edges of the
hexagon. There are only eleven such 4-sets, namely, $x135$,
$x246$, $1235$, $2346$, $1345$, $2456$, $1356$, $1246$, $1245$,
$2356$, $1346$. Since at most two of the four sets $x135$, $1235$,
$1345$, $1356$ can be facets, and at most two of the four sets
$x246$, $2346$, $2456$, $1246$ can be facets, we have no way to
choose nine of these eleven sets as facets of $M^3_9$. So,
$M^3_9$ has no edge of degree 6.

Thus all the edges of $M^3_9$ have degree 3, 4 or 5. For $3 \leq
i \leq 5$, let $\varepsilon_i$ be the number of edges of degree
$i$. Since the total number of edges is ${9 \choose 2} = 36$, we
have
\begin{equation} \label{eq2}
\varepsilon_3 + \varepsilon_4 + \varepsilon_5 = 36.
\end{equation}
Also, counting in two ways the ordered pairs $(e, \sigma)$, where
$e$ is an edge in a facet $\sigma$, we get
\begin{equation} \label{eq3}
3\varepsilon_3 + 4\varepsilon_4 + 5\varepsilon_5 = 27 \times {4
\choose 2} = 162.
\end{equation}
Since any two facets intersect, Equation (\ref{28/29}) shows that
${\displaystyle \sum_{e\subset \sigma} \deg(e)}=28$ for each
facet $\sigma$. Since the only permissible edge-degrees are 3, 4
and 5, it follows that there are only two types of facets. A
facet of type 1 contains one edge of degree 3 (and five of degree
5) while a facet of type 2 contains two edges of degree 4 (and
four of degree 5). Counting in two ways pairs $(e, \sigma)$ with
$e\subset \sigma$, where (i) $e$ is an edge of degree 3 and
$\sigma$ is a facet (of type 1) and (ii) $e$ is an edge of degree
4 and $\sigma$ is a facet (of type 2), we see that there are
$3\varepsilon_3$ facets of type 1 and $2\varepsilon_4$ facets of
type 2. Since the total number of facets is 27, we get
\begin{equation} \label{eq4}
3\varepsilon_3 + 4\varepsilon_4 = 27.
\end{equation}
Solving Equations (\ref{eq2}), (\ref{eq3}) and (\ref{eq4}), we
obtain $\varepsilon_3 = 9$, $\varepsilon_4 = 0$, $\varepsilon_5 =
27$. Thus all the edges have degree 3 or 5. Therefore, for any
vertex $x$, ${\rm lk}(x)$ is an $S^{\,2}_8$ all whose vertices
have degree 3 or 5. Since an $S^{\,2}_8$ has 18 edges, its vertex
degrees add up to 36. So, ${\rm lk}(x)$ has exactly two vertices
of degree 3. That is, each vertex $x$ of $M^3_9$ is in exactly two
edges of degree 3. Hence, by Lemma \ref{l4.3}, there is an edge of
degree $\geq 6$, a contradiction. This proves the lemma. \hfill
$\Box$

\bigskip

Let's say that a vertex $x$ of $M^3_9$ is {\em good} if there is
a partition of $V(M^3_9) \setminus\{x\}$ into two facets. By,
Lemma \ref{l4.4}, there is at least one good vertex. Next we
prove.

\begin{lemma}$\!\!${\bf .} \label{l4.5}
The link of any good vertex in $M^3_9$ is isomorphic to the
$2$-sphere ${\cal S}$ given below.
\end{lemma}

\noindent {\bf Proof.} Let $v$ be a good vertex. Let $1235$ and
$4678$ be two disjoint facets not containing $v$. Let $L = {\rm
lk}(v)$. In view of Lemma \ref{l4.2}, we may assume that the
induced subcomplex of $L$ on $1235$ and $4678$ are
$S^{\,1}_3(\{1, 2, 3\}) \cup \{5\}$ and $S^{\,1}_3(\{6, 7, 8\})
\cup \{4\}$ respectively. Hence $V({\rm lk}_{L}(5)) \subseteq
\{4, 6, 7, 8 \}$ and $V({\rm lk}_{L}(4)) \subseteq \{1, 2, 3,
5\}$. It follows that no triangle of $L$ contains $45$, so that
$45$ is not an edge of $L$. Therefore, ${\rm lk}_{L}(5) =
S^{\,1}_3(\{6, 7, 8\})$ and ${\rm lk}_{L}(4) = S^{\,1}_3(\{1, 2,
3\})$. Since each vertex of $L$ is adjacent in $L$ with 4 or 5,
and no two degree 3 vertices are adjacent in an $S^{\,2}_8$, it
follows that 4 and 5 are the only two vertices of degree 3 in $L$.
Note that the partition $\{1235, 4678\}$ of
$V(M^3_9)\setminus\{v\}$ into two facets of $M^3_9$ is uniquely
recovered from $L$ as the pair of stars of the degree 3 vertices
of $L$. (Star of a vertex $u$ in a simplicial complex $K$ is the
join $\{u\} \ast {\rm lk}_K(u)$.) Thus, there is a natural
bijection between good vertices and pairs of disjoint facets.

Collapsing the two degree 3 vertices 4 and 5 in $L$, we obtain an
$S^{\,2}_6$ with two disjoint triangle $123$ and $678$. Therefore,
$L$ is obtained from an $S^{\,2}_6$ by starring a vertex in each
of two disjoint triangles. We know (see Proposition \ref{s27})
that there are exactly two different $S^{\,2}_6$, namely, ${\cal
S}_3$ and ${\cal S}_4$. Observe that each of these two
$S^{\,2}_6$ has a unique pair of disjoint triangles, up to
automorphisms of the $S^{\,2}_6$. Thus $L$ is isomorphic to
${\cal S}$ or ${\cal T}$ below.

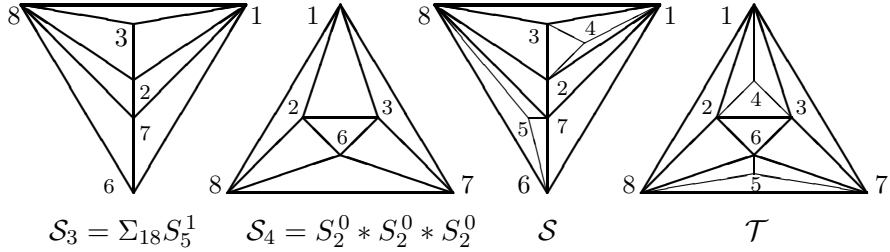
\setlength{\unitlength}{2.5mm}

\begin{picture}(47,14.5)(0,0)

\thicklines

\put(7,3){\line(-3,5){6}} \put(7,3){\line(0,1){9}}
\put(7,3){\line(3,5){6}} \put(1,13){\line(1,0){12}}

\put(7,7){\line(-1,1){6}} \put(7,7){\line(1,1){6}}
\put(7,9){\line(-3,2){6}} \put(7,9){\line(3,2){6}}
\put(7,12){\line(-6,1){6}} \put(7,12){\line(6,1){6}}

\put(0.3,12){\small 8} \put(13.2,12){\small 1}
\put(6,10.8){\small 3} \put(7.3,8){\scriptsize 2}
\put(7.3,5.8){\scriptsize 7} \put(5.4,3){\scriptsize 6}

\put(2.5,0.5){\mbox{${\cal S}_3 = \Sigma_{18}S^{\,1}_5$}}


\thicklines

\put(12,3){\line(3,5){6}} \put(12,3){\line(1,1){4}}
 \put(12,3){\line(3,1){6}} \put(12,3){\line(1,0){12}}

\put(24,3){\line(-3,5){6}} \put(24,3){\line(-1,1){4}}
\put(24,3){\line(-3,1){6}} \put(18,5){\line(1,1){2}}
\put(18,5){\line(-1,1){2}} \put(16,7){\line(1,0){4}}
\put(16,7){\line(1,3){2}} \put(20,7){\line(-1,3){2}}

\put(11,3){\small 8} \put(24.4,3){\small 7}
 \put(17.8,5.6){\scriptsize 6} \put(15.2,7.2){\scriptsize 2}
 \put(20.2,7.2){\scriptsize 3} \put(16.2,12){\small 1}

\put(13,0.5){\mbox{${\cal S}_4 = S^{\,0}_2\ast S^{\,0}_2 \ast
S^{\,0}_2$}}


\thicklines

\put(29,3){\line(-3,5){6}} \put(29,3){\line(0,1){9}}
\put(29,3){\line(3,5){6}} \put(23,13){\line(1,0){12}}

\put(29,7){\line(-1,1){6}} \put(29,7){\line(1,1){6}}
\put(29,9){\line(-3,2){6}} \put(29,9){\line(3,2){6}}
\put(29,12){\line(-6,1){6}} \put(29,12){\line(6,1){6}}

\thinlines

\put(29,3){\line(-1,4){1}} \put(28,7){\line(1,0){1}}
\put(28,7){\line(-5,6){5}} \put(29,9){\line(1,1){2}}
\put(31,11){\line(2,1){4}} \put(31,11){\line(-2,1){2}}

\put(22.1,12){\small 8} \put(35.2,12){\small 1}
\put(28,11){\scriptsize 3}  \put(31,11.4){\scriptsize 4}

\put(29.5,8.2){\scriptsize 2} \put(29.5,6.2){\scriptsize 7}
\put(27.4,6){\scriptsize 5} \put(27.4,3){\small 6}

\put(28.5,0.5){\mbox{${\cal S}$}}


\thicklines

\put(34,3){\line(3,5){6}} \put(34,3){\line(1,1){4}}
 \put(34,3){\line(3,1){6}} \put(34,3){\line(1,0){12}}

\put(46,3){\line(-3,5){6}} \put(46,3){\line(-1,1){4}}
\put(46,3){\line(-3,1){6}} \put(40,5){\line(1,1){2}}
\put(40,5){\line(-1,1){2}} \put(38,7){\line(1,0){4}}
\put(38,7){\line(1,3){2}} \put(42,7){\line(-1,3){2}}

\thinlines

\put(34,3){\line(6,1){6}} \put(46,3){\line(-6,1){6}}
\put(40,4){\line(0,1){1}} \put(38,7){\line(1,1){2}}
\put(42,7){\line(-1,1){2}} \put(40,9){\line(0,1){4}}

\put(33,3){\small 8} \put(39.8,3.1){\scriptsize 5}
\put(46.4,3){\small 7} \put(39.8,5.6){\scriptsize 6}
\put(39.8,7.6){\scriptsize 4} \put(37.2,7.2){\scriptsize 2}
\put(42.2,7.2){\scriptsize 3} \put(38.2,12){\small 1}

\put(39.5,0.5){\mbox{${\cal T}$}}

\end{picture}

If possible, let $L = {\cal T}$. We claim that the facet of
$M^3_9$ (other than $v124$) containing $124$ must be $1245$.
Indeed, it can not be $1234$ since then there would be three
facets of $M^3_9$ disjoint from it, contradicting Lemma
\ref{l4.2}. Also, it can not be $124i$ for $6\leq i\leq 8$ (since
the induced subcomplex on its complement contradicts Lemma
\ref{l4.2}). Thus, $1245$ is a facet of $M^3_9$. Similarly, we
get six facets $1245, 1345, 2345, 4567, 4568, 4578$ of $M^3_9$.
Then ${\rm lk}_{M^3_9}(45)$ is the disjoint union of two circles.
This is not possible since $M^3_9$ is a manifold. \hfill $\Box$

\begin{lemma}$\!\!${\bf .} \label{l4.6}
If $M^3_9$ is a $9$-vertex non-sphere neighbourly combinatorial
$3$-manifold then $M^3_9 = K^3_9$.
\end{lemma}

\noindent {\bf Proof.} By Lemma \ref{l4.4}, there is a good
vertex, say 9 is a good vertex. By Lemma \ref{l4.5}, the link of
9 is isomorphic to ${\cal S}$. Assume that the link of 9 is
${\cal S}$. The facet (other than $2349$) containing $234$ must
be $2346$ (since for every vertex $x\neq 6, 9$, there are two
facets through 9 disjoint from $234x$). So, $2346$ and $5789$ are
disjoint facets. Then 1 is a good vertex. Similarly, $3567\in
M^3_9$ and $8$ is a good vertex. A similar argument shows that
the facet (other than $1349$) through $134$ must be $1345$ (since
the induced subcomplex of ${\rm lk}(2)$ on $5689$ is not an
$S^{\,1}_3$ together with an isolated vertex, by Lemma
\ref{l4.2}, $1347$ can not be a facet). Similarly, the facet
(other than $5689$) through $568$ must be $4568$.

Now, consider the links of 1 and 8. By Lemma \ref{l4.5}, both are
isomorphic to ${\cal S}$. Note that the only non-neighbour in
${\cal S}$ of a vertex of degree 6 has degree 3. Since ${\rm
lk}(89) =$ \setlength{\unitlength}{2mm}
\begin{picture}(5,2.5)(4,-0.5)

\thicklines

\put(3,0){\line(2,1){2}} \put(3,0){\line(2,-1){2}}
\put(5,1){\line(1,0){1.6}} \put(5,-1){\line(1,0){1.6}}
\put(6.6,1){\line(2,-1){2}} \put(6.6,-1){\line(2,1){2}}

\put(2.8,0.7){\scriptsize 1} \put(4.3,1.3){\scriptsize 3}
\put(6.5,1.3){\scriptsize 2} \put(8.2,0.7){\scriptsize 7}
\put(6.5,-0.6){\scriptsize 5} \put(4.7,-0.6){\scriptsize 6}
\end{picture} and
${\rm lk}(19) =$ \setlength{\unitlength}{2mm}
\begin{picture}(5,2.5)(4,-0.5)

\thicklines

\put(3,0){\line(2,1){2}} \put(3,0){\line(2,-1){2}}
\put(5,1){\line(1,0){1.6}} \put(5,-1){\line(1,0){1.6}}
\put(6.6,1){\line(2,-1){2}} \put(6.6,-1){\line(2,1){2}}

\put(2.8,0.7){\scriptsize 2} \put(4.3,1.3){\scriptsize 4}
\put(6.5,1.3){\scriptsize 3} \put(8.2,0.7){\scriptsize 8}
\put(6.5,-0.6){\scriptsize 6} \put(4.7,-0.6){\scriptsize 7}
\end{picture},
it follows that $\deg(48) = 3 = \deg(15)$. Since $3567$ and $1249$
are disjoint facets not containing 8 and since $189$, $289\in
M^3_9$, Lemma \ref{l4.2} implies that $128\in M^3_9$ and $148
\not\in M^3_9$. Since ${\rm lk}(8)$ is a copy of ${\cal S}$ and
$12$ is an edge and $14$ is a non-edge in ${\rm lk}(8)$, it
follows that $123$ is a triangle of ${\rm lk}(8)$ and hence
$\deg_{{\rm lk}(8)}(3) = 3$. Since the two degree 3 vertices are
non-adjacent in the edge-graph of ${\cal S}$ and $\deg_{{\rm
lk}(8)}(3) = 3 = \deg_{{\rm lk}(8)}(4)$, it follows that ${\rm
lk}_{{\rm lk}(8)}(4) = S^{\,1}_3(5, 6, 7)$. Similarly, since
$2346$ and $5789$ are disjoint facets not containing 1 and $189$,
$179\in M^3_9$, it follows that $178\in M^3_9$, $158\not\in
M^3_9$ and hence ${\rm lk}_{{\rm lk}(1)}(5) = S^{\,1}_3(2, 3,
4)$. Then, from the links of 1 and 8 we get facets $1235$,
$1245$, $1238$, $1278$, $1678$, $4578$, $4678$.

Now, trying to complete the links of 2 and 7, we get facets
$2356$, $2456$, $3457$, $3467$. This implies that $K^{\,3}_9$ is a
subcomplex of $M^{\,3}_9$. Since both are 3-pseudomanifolds,
$M^{\,3}_9 = K^{\,3}_9$. \hfill $\Box$

\bigskip

\noindent {\bf Proof of Theorem \ref{t2}.} Observe that the
degree 3 edges in $K^{\,3}_9$ are $15$, $59$, $94$, $48$, $83$,
$37$, $72$, $26$, $61$ and the automorphism group $D_{18}$ of
$K^{\,3}_9$ acts transitively on these nine edges. But, none of
these nine edges are removable. (Since ${\rm lk}_{K^{\,3}_9}(15) =
S^{\,1}_3(2, 3, 4)$ and $234$ is a face in $K^{\,3}_9$, $15$ is
not removable.) So, there is no bistellar 2-move on $K^{\,3}_9$.

Now, let $N^{\,3}_9$ be a 9-vertex non-sphere combinatorial
3-manifold. If $N^{\,3}_9$ is not neighbourly then, by Theorem
\ref{t1}, there is a 9-vertex neighbourly 3-manifold $M^{\,3}_9$
obtainable from $N^{\,3}_9$ by a sequence of bistellar 1-moves.
Since $N^{\,3}_9$ is non-sphere, so is $M^{\,3}_9$. Therefore, by
Lemma \ref{l4.6}, $M^{\,3}_9 = K^{\,3}_9$. Thus, $N^{\,3}_9$ is
obtainable from $K^{\,3}_9$ by a sequence of bistellar 2-moves.
But, we just observe that $K^{\,3}_9$ does not admit any
bistellar 2-move. A contradiction. So, $N^{\,3}_9$ is
neighbourly. Now, by Lemma \ref{l4.6}, $N^{\,3}_9 = K^{\,3}_9$.
\hfill $\Box$

\bigskip

\noindent {\bf Acknowledgement\,:} B. Datta thanks Department of
Science \& Technology (DST), India for financial support through a
project (Grant: SR/S4/MS-272/05).

\end{document}